\documentclass[leqno,twoside]{article}
\usepackage{amsmath,amssymb,amsthm,extarrows}
\usepackage[margin=2.5cm]{geometry}

\usepackage{titlesec,hyperref}
\usepackage{fancyhdr}

\fancypagestyle{plain}
{
\fancyhf{}

}

\titleformat{\subsection}{\it}{\thesubsection.\enspace}{1pt}{}

\newtheorem{theo}{Theorem}[section]
\newtheorem{lemm}[theo]{Lemma}
\newtheorem{defi}[theo]{Definition}

\newtheorem{rema}[theo]{Remark}
\numberwithin{equation}{section}

\allowdisplaybreaks

\begin{document}
\title{Gevrey Regularity for Solutions of the Non-Cutoff Boltzmann Equation: Spatially Inhomogeneous Case
\hspace{-4mm}
}

\author{Teng-Fei Zhang$^1$
\quad Zhaoyang Yin$^2$ \\[10pt]
Department of Mathematics, Sun Yat-sen University,\\
510275, Guangzhou, P. R. China.}
\footnotetext[1]{Corresponding author. Email: {\it fgeyirui@163.com}}
\footnotetext[2]{Email: \it mcsyzy@mail.sysu.com.cn}
\date{}
\maketitle
\hrule

\begin{abstract}

In this paper we consider the non-cutoff Boltzmann equation in spatially inhomogeneous case. We prove the propagation of Gevrey regularity for the so-called smooth Maxwellian decay solutions to the Cauchy problem of spatially inhomogeneous Boltzmann equation, and obtain Gevrey regularity of order $1/s$ in the velocity variable $v$ and order $1$ in the space variable $x$. The strategy relies on our recent results for spatially homogeneous case (J. Differential Equations 253(4) (2012), 1172-1190. DOI: 10.1016/j.jde.2012.04.023). Rather, we need much more intricate analysis additionally in order to handle with the coupling of the double variables. Combining with the previous result mentioned above, it gives a characterization of the Gevrey regularity of the particular kind of solutions to the non-cutoff Boltzmann.

\vspace*{5pt}
\noindent {\it 2000 Mathematics Subject Classification}: 35A05, 35B65, 35D10, 35H20, 76P05, 82C40.

\vspace*{5pt}
\noindent{\it Keywords}: Boltzmann equation; Spatially inhomogeneous; Non-cutoff; Gevrey regularity.
\end{abstract}

\vspace*{10pt}
\tableofcontents

\section{Introduction}

\subsection{The Boltzmann equation}

In this paper we consider the Cauchy problem of the spatially inhomogeneous Boltzmann equation without angular cutoff. It reads, with a $T >0$, as the following equation,
\begin{align}\label{BE}
  \left\{
    \begin{array}{l}
      f_t(t,x,v)+v\cdot \nabla_xf(t,x,v)=Q(f,f)(v),~~ t\in (0,T], \\
      f(0,x,v)=f_0(x,v),
    \end{array}
  \right.
\end{align}
for the density distribution function of particles $f=f(t,x,v)$, which are located around position $x \in \mathbb{T}^3$ with velocity $v\in \mathbb{R}^3$ at time $t \geq 0$. The right-hand side of the above equation is the so-called Boltzmann bilinear collision operator acting only on the velocity variable $v$:
\begin{align*}
  Q(g, f)=\int_{\mathbb{R}^3}\int_{\mathbb S^{2}}B\left({v-v_*},\sigma \right)
    \left\{ g'_* f'-g_*f \right\} d\sigma dv_*.
\end{align*}

Above, we use the standard shorthand $f=f(t,x,v)$, $f_*=f(t,x,v_*) $, $f'=f(t,x,v') $, $f'_*=f(t,x,v'_*) $. The relations between the post- and pre-collisional velocities are described by the $\sigma$-representation, that is, for $\sigma \in \mathbb S^2$,
$$
v'=\frac{v+v_*}{2}+\frac{|v-v_*|}{2}\sigma,~~ v'_*
=\frac{v+v_*}{2}-\frac{|v-v_*|}{2}\sigma.
$$
Note that the collision process satisfies the conservation of momentum and kinetic energy, i.e.
$$
v+v_*=v'+v'_*,\qquad  |v|^2+|v_*|^2=|v'|^2+|v'_*|^2.
$$

The collision cross section $B(z, \sigma)$ is a given non-negative function depending only on the interaction law between particles. From a mathematical viewpoint, that is to say, $B(z, \sigma)$ depends only on the relative velocity $|z|=|v-v_*|$ and the deviation angle $\theta$ defined through the scalar product $\cos \theta=\frac{z}{|z|} \cdot \sigma.$

Without loss of generality, the cross section $B$  is assumed to be of the form:
$$
B(v-v_*, \cos \theta)=\Phi (|v-v_*|) b(\cos \theta),~~
\cos \theta=\frac{v-v_*}{|v-v_*|} \cdot \sigma , ~~
0\leq\theta\leq\frac{\pi}{2},
$$
where the kinetic factor $\Phi$ is given by
$$
\Phi (|v-v_*|) = |v-v_*|^{\gamma},
$$
and the angular part $b$, with a singularity, satisfies,
$$
\sin \theta b(\cos \theta) \sim \theta^{-1-2s}, \ \ \mbox{as} \ \ \theta\rightarrow 0+,
$$
for some $0< s <1$.

We remark that if the inter-molecule potential is given by the inverse-power law $U(\rho) = \rho ^{-(p-1)}$ (where $p>2$), it holds $ \gamma = \frac{p-5}{p-1} $, $ s=\frac{1}{p-1}.$ Generally, the cases $\gamma >0$, $\gamma =0$, and $\gamma <0$ correspond to so-called hard, Maxwellian, and soft potential respectively. And the cases $0<s<1/2$, $1/2 \leq s<1$ correspond to so-called mild singularity and strong singularity respectively.

\subsection{Review of non-cutoff theory in Gevrey spaces}

We begin with a brief review for the non-cutoff theory of the Boltzmann equation. We refer to Villani's review book \cite{Villani} for the physical background and the mathematical theories of the Boltzmann equation. Furthermore, in the non-cutoff setting, Alexandre gave more details in \cite{Alex-review}.

Our discussion is based on the following definition of Gevrey spaces $G^s(\Omega)$ on an open subset $\Omega \subseteqq \mathbb{R}^3$ (see \cite{Luigi}, for instance):
\begin{defi}\label{Def1}
For $0<s<+\infty$, we say that $ f \in G^s(\Omega) $, if $f \in C^\infty(\Omega)$, and there exist $C>0,~ N_0>0$ such that
$$ \|\partial^\alpha f\|_{L^2(\Omega)} \leq C^{|\alpha|+1} {\{\alpha! \}^s},\quad
\forall \alpha \in \mathbb{N}^3,~ |\alpha| \geq N_0.
$$
\end{defi}

Note that the Gevrey scale measures regularity between analytic and $C^\infty$. More precisely, when $s = 1$, it is usual analytic function. If $s > 1$, it is Gevrey class function. And for $0 < s < 1$, it is called ultra-analytic function.

For the Cauchy problem of the Boltzmann equation in Gevrey classes, Ukai showed, in \cite{Ukai-84} in 1984, that there exists a unique local solution for both spatially homogeneous and inhomogeneous cases, with the assumption on the cross section:
\begin{align*}
&\big| B(|z|,\cos \theta) \big| \leq  K(1+|z|^{-\gamma'}+|z|^\gamma) \theta^{-n+1-2s},
      \quad n  \textrm{ is dimensionality},\\
&(0\leq \gamma' < n,~ 0\leq \gamma <2,~ 0\leq s<1/2, ~\gamma +6s<2 ).
\end{align*}

In particular, for the spatially inhomogeneous case, by introducing the norm of Gevrey space
$$
\|f\|_{\delta,\rho_1,\nu_1,\rho_2,\nu_2} = \sum_{\alpha,\beta} \frac{\rho_1^{|\alpha|} \rho_2^{|\beta|}}
          {\{\alpha!\}^{\nu_1} \{\beta!\}^{\nu_2}}
\|e^{\delta \langle v \rangle^2} \partial_x^{\alpha} \partial_v^{\beta} f \|_{L^\infty(\mathbb{R}^n_x \times \mathbb{R}^n_v)},
$$
Ukai proved that, under some assumptions for $\nu$ and the initial datum $f_0(x,v)$, the Cauchy problem (\ref{BE}) has a unique solution $f(t,x,v)$ for $t\in (0,T]$.

On the other hand, Desvillettes established in \cite{Desvillettes2} the $C^\infty$ smoothing effect for solutions of Cauchy problem in spatially homogeneous case, and conjectured Gevrey smoothing effect. He also proved, without any assumptions on the decay at infinity in $v$ variables, the propagation of Gevrey regularity for solutions (see \cite{Desvillettes}).

In 2009 Morimoto et al. considered in \cite{Mori-Ukai-Xu-Yang} the Gevrey regularity for the linearized Boltzmann equation around the absolute Maxwellian distribution, by virtue of the following mollifier:
$$
G_\delta (t,D_v)=\frac{e^{t \langle D_v \rangle^{1/\nu}}}
                      {1+\delta e^{t \langle D_v \rangle^{1/\nu}}},\quad 0< \delta <1.
$$
We remark that the same operator was used in many related models such as the Fokker-Planck equation, the Kac's equation, the Landau equation, and so on.

In the mild singularity setting $0<s<1/2$, Huo et al. proved in \cite{Huo} that any weak solution $f(t,v)$ to the Cauchy problem (\ref{BE}) satisfying the natural boundedness on mass, energy and entropy, namely,
\begin{align}\label{natural bound}
\int_{\mathbb{R}^n} f(v)[1+|v|^2+\log(1+f(v))]dv < +\infty,
\end{align}
belongs to $H^{+\infty}(\mathbb{R}^n)$ for any $0<t \leq T$, and moreover,
\begin{align}\label{smooth solution}
f \in L^\infty \big([t_0,T];H^{+\infty}(\mathbb{R}^n)\big),
\end{align}
for any $T>0$ and $t_0 \in (0,T)$.

In paper \cite{FiveGroup-regulariz}, the five authors proved the smoothing effect on the solution with weight. More percisely, if the non-negative $f$ belongs to $\mathcal{H}_l^5 \Big( (t_1,t_2) \times \Omega \times \mathbb{R}^3_v \Big)$, solves the Boltzmann equation (\ref{BE}) in the above domain in the classic sense, and meantime satisfies the non-vacuum condition $\|f(t,x,v)\|_{L^1(\mathbb{R}^3_v)} >0$, then it follows that,
$$
f\in \mathcal{H}_l^\infty \Big( (t_1,t_2) \times \Omega \times \mathbb{R}^3_v \Big),
$$
and hence it holds,
$$
f\in C^\infty \Big( (t_1,t_2) \times \Omega; \mathcal{S}(\mathbb{R}^3_v) \Big).
$$

Therein the five authors also considered a kind of solution having the Maxwellian decay, based on which we introduce the following definition:
\begin{defi}\label{Def2}
We say that $f(t,x,v)$ is a smooth Maxwellian decay solution to the Cauchy problem (\ref{BE}) if
\begin{align*}
  \left\{
     \begin{array}{l}
       f \geq 0,~ \not \equiv 0,\\
       \exists ~ \delta_0 >0 \textrm{ such that }
          e^{\delta_0 \langle v \rangle ^2}
          f \in L^\infty \left( [0,T];~H^{+\infty}(\mathbb{R}^3_x \times \mathbb{R}^3_v) \right).
     \end{array}
  \right.
\end{align*}
\end{defi}

Note that the Theorem 1.2 of \cite{FiveGroup-regulariz} shows the uniqueness of the smooth Maxwellian decay solution to the Cauchy problem (\ref{BE}).

In 2010 Morimoto-Ukai considered the Gevrey regularity of $C^\infty$ solutions with the Maxwellian decay to the Cauchy problem of spatially homogeneous Boltzmann equation (see \cite{Mori-Ukai}). Motivated by their results, we studied this problem in \cite{Zhang-Yin} in a more general framework. More precisely, we considered the general kinetic factor $ \Phi (|v|) = |v|^{\gamma} $ instead of the moderate form $ \langle v \rangle ^{\gamma}=( 1+|v|^2)^{\gamma/2} $ in \cite{Mori-Ukai}, and a wider range of the parameter of $\gamma$ such that $\gamma +2s \in (-1,1)$ which applies for both hard potential and soft potential.

In the present paper, we study still in the general framework $\Phi (|v|) = |v|^{\gamma}$ with $\gamma +2s \in (-1,1)$ in the mild singularity assumption $0<s<1/2$. Beyond that, we focus here on the spatially inhomogeneous case, which is much more complicate than homogeneous case (because, in the spatially inhomogeneous case, the interaction between the kinetic part and nonlinear collision part is very complicate). The full, spatially inhomogeneous model is more closely related to the real physical setting, thus, is more meaningful and interested, and in particular is a cornerstone of statistical physics.

Since the estimates obtained in \cite{Zhang-Yin} can carry over to the inhomogeneous framework, it seems easy to begin our justification. We can handle with the space variables $x$ by virtue of one more integrations. However, it is difficult to clarify distinctly the process when taking double supremum on space variables $x$ and velocity variables $v$, which will be related to a more intricate technique, as we will see later. We aim in this work at expressing the whole process explicitly. Combining with the previous result in \cite{Zhang-Yin}, we can get a characterization of the Gevrey regularity of smooth Maxwellian decay solutions to the non-cutoff Boltzmann. We wish these results will be useful for the forthcoming research.

\subsection{Main results}

When considering the Gevrey regularity we may assume $t_0=0$ in the above statement by translation.

Now we are in a position to state our main result of propagation of Gevrey regularity, as follows:
\begin{theo}\label{Propagation of Gevrey}
Let $\nu_1 \geq 1,~\nu_2 >1 $(which are independent of s) and assume that $ 0<s<1/2 $, $ -1<\gamma+2s<1$. Let $f(t,x,v)$ be a smooth Maxwellian decay solution to the Cauchy problem (\ref{BE}). If there exist $\rho'$, $\delta'$ such that
\begin{align}\label{initial assumption}
\sup_{\alpha,\beta} \frac{\rho'^{|\alpha|+|\beta|} \|e^{\delta' \langle v \rangle ^2}\partial_x^\alpha \partial^{\beta}_{v} f(0)\|_{L^2_{x,v}}}{\{\alpha!\}^{\nu_1} \{\beta!\}^{\nu_2}} < +\infty,
\end{align}
then there exist $ \rho >0 $ and $ \delta,\kappa >0 $ with $ \delta >\kappa T $ such that
\begin{align}\label{propag}
\sup_{t\in (0,T]}\sup_{\alpha,\beta} \frac{\rho^{|\alpha|+|\beta|} \|e^{(\delta-\kappa t) \langle v \rangle ^2}\partial_x^\alpha \partial^{\beta}_{v} f(t)\|_{L^2_{x,v}}}{\{\alpha!\}^{\nu_1} \{\beta!\}^{\nu_2}} < +\infty.
\end{align}
\end{theo}

\begin{rema}
It should be noted that the above theorem is similar as Theorem 1.3 in \cite{Zhang-Yin}, but here we consider the spatially inhomogeneous case.
\end{rema}

By using the arguments in Section 4 below, we obtain the Gevrey smoothing effect of order $1/s$ in variable $v$ and order $1$ in variable $x$ as follows:
\begin{theo}\label{Gevrey Regularity}
Assume that $0<s<1/2,\ -1<\gamma+2s<1$ and $\nu_1=1,~\nu_2 =1/s$. Let $f(t,x,v)$ be a smooth Maxwellian decay solution to the Cauchy problem (\ref{BE}), and further, $f(0,x,v)$ be analytic with respect to space variable $x$, then for any $t_0 \in (0,T)$, there exist $\rho >0 $ and $\delta,\kappa >0 $ with $ \delta >\kappa T $ such that
\begin{align}
\sup_{t\in [t_0,T]}\sup_{\alpha,\beta} \frac{\rho^{|\alpha|+|\beta|} \|e^{(\delta-\kappa t) \langle v \rangle ^2}\partial^{\alpha}_{x} \partial_v^\beta f(t)\|_{L^2_{x,v}}}{\{\alpha!\}^{\nu_1} \{\beta!\}^{\nu_2}} < +\infty.
\end{align}
\end{theo}

\subsection{The structure of the paper}

The remainder of the paper proceeds as follows. In the next section we give some preliminaries and a key lemma, by which we can complete immediately the proof of Theorem \ref{Propagation of Gevrey}. Section 3 is devoted to the proof of the key lemma. Section 4 is arranged for the justification of Theorem \ref{Gevrey Regularity} about the Gevrey smoothing effect.

\section{Preliminaries and the key lemma}

\subsection{Preliminaries}

We give some notations and fundamental facts here. (see \cite{Mori-Ukai,Zhang-Yin} for details).

Denote $\partial^\alpha_\beta f = \partial_x^\alpha \partial_v^\beta f$ throughout this paper.
Let $l,r \in \mathbb{Z}_+$ whose values will be chosen later. For $\delta,\rho>0$ we define:
\begin{align}\label{double norm}
   \| f \|_{\delta,l,\rho,r,\alpha,\beta} \xlongequal{\!\!def\!\!}
   \frac{\rho^{|\alpha|+|\beta|} \| \langle v \rangle^l e^{\delta \langle v \rangle^2}
   \partial^\alpha_\beta f\|_{L^2_x(\mathbb{T}^3) L^2_v(\mathbb{R}^3)} }
   {\{ (\alpha-r)!\}^{\nu_1} \{ (\beta-r)!\}^{\nu_2}},
\end{align}
where $\alpha,~\beta$ are multi-index of $\mathbb{Z}^3$, i.e., $ \alpha =(\alpha^1,\alpha^2,\alpha^3) \in \mathbb{Z}^3_+,~\beta =(\beta^1,\beta^2,\beta^3) \in \mathbb{Z}^3_+$. We denote
$$ (\alpha - r)!=(\alpha^1 - r)! (\alpha^2 - r)! (\alpha^3 - r)!,$$
and
$$ C_\alpha^{\alpha_1}=\frac{\alpha!}{\alpha_1! \alpha_2!}.$$

We introduce the definition:
\begin{align}
  \| f \|_{l,\rho,r,N} \xlongequal{\!\!def\!\!}
  \sup_{8r \leq |\alpha|+|\beta| \leq N}
  \| f \|_{\delta-\kappa t,l,\rho,r,\alpha,\beta}
\end{align}
with fixed $\delta,~\kappa >0$ satisfying $\delta > \kappa T$. Here $N$ is a fixed large integer.

Analogous to the argument of \cite{Mori-Ukai}, we can obtain that, for $h>1$,
\begin{align}
\|f\|_{l,\rho (1+h)^{-\nu_1-\nu_2},r,N}(t) \leq \Big\{\frac{(r!)^3}{h^r}\Big\}^{\nu_1+\nu_2} \|f\|_{l,\rho,0,N}(t).
\end{align}

Setting $\rho =\rho'$ and taking a large enough $ h $, then it follows from the assumption (\ref{initial assumption}) that, we may let $ \|f\|_{l,\rho'(1+h)^{-\nu_1-\nu_2},r,N}(0) $ be as small as possible, where $\delta$ can be chosen any positive less than $\delta' > 0$ in (\ref{initial assumption}).

Therefore, in order to prove the result (\ref{propag}), it suffices to prove that,
\begin{align}\label{Supremum}
\sup_{t \in (0,T]} \|f\|_{l,\rho,r,N}(t) < \infty,
\end{align}
under the assumption that $ \|f\|_{l,\rho,r,N}(0) $ is sufficiently small. Above, $\rho =\rho' (1+h)^{-\nu_1-\nu_2}$.

\subsection{The key lemma}

Now we give the key lemma, which will play an important role in the following sections.
\begin{lemm}\label{Key lemma}
If $ l \geq 4$ and $r>1+\nu_2/(\nu_2-1)$, then for any $\alpha,~\beta$ satisfying $ 8r \leq |\alpha|+|\beta| \leq N$ we have
\begin{align}\label{basic inequa}
&\|f(t)\|^2_{\delta -\kappa t,l,\rho,r,\alpha,\beta}
+ 2\kappa \int_0^t  \|f(\tau)\|^2_{\delta-\kappa \tau,l+1,\rho,r,\alpha,\beta}d\tau \\
\leq
&\|f(0)\|^2_{\delta,l,\rho,r,\alpha,\beta}
       + C_\kappa \int^t_0 \Big( \|f\|^2_{l,\rho,r,N}(\tau)
                             + \|f\|^{2(1+\eta)/\eta}_{l,\rho,r,N}(\tau) \Big)d\tau \nonumber\\
&+ \frac{\kappa}{10}\sup_{8r \leq |\alpha|+|\beta| \leq N}
                         \int^t_0\! \|f(\tau)\|^2_{\delta-\kappa \tau,l+1,\rho,r,\alpha,\beta}d\tau, \nonumber
\end{align}
where $\eta \in (0,1)$.
\end{lemm}

Following along the same lines as that of Section 2 in \cite{Mori-Ukai}, we can prove Theorem \ref{Propagation of Gevrey}, and so omit it.

The proof of this lemma will be given in the next section.

\section{The proof of the key lemma}

\subsection{Rewrite the equation}

Applying $\partial^\alpha_\beta$ to Eq.(\ref{BE}), we have, for $|\beta|=0$,
\begin{align}
  \partial_t (\partial^\alpha_\beta f) + v \cdot \nabla_x (\partial^\alpha_\beta f)
  = \partial^\alpha_\beta Q(f,f),
\end{align}
and for $|\beta| \geq 1$,
\begin{align}
  \partial_t (\partial^\alpha_\beta f) + v \cdot \nabla_x (\partial^\alpha_\beta f)
  = -\sum_j \beta^j (\partial^{\alpha+e_j}_{\beta-e_j} f) + \partial^\alpha_\beta \,Q(f,f),
\end{align}
where $j=1,2,3$, and $\beta=(\beta^1,\beta^2,\beta^3),~e_1=(1,0,0),e_2=(0,1,0),e_3=(0,0,1)$. Note that if $\beta^j-1$, the $j$-th component of $\beta-e_j$, equals to $-1$, then $\beta^j$ will be zero, the above equation will still hold true. %We will consider only the case $|\beta| \geq 1$ because the estimates for the case $|\beta|=0$ are similar and easier.

Let $\mu =\mu_{\delta,\kappa}(t) =e^{-(\delta-\kappa t) \langle v
\rangle ^2}  $ with $\delta > \kappa T$. Multiplying by $\mu^{-1}$
both sides of the above equation, we obtain
\begin{align}\label{eq1}
\partial_t (\mu^{-1} \partial^\alpha_\beta f) + v\cdot \nabla_x (\mu^{-1} \partial^\alpha_\beta f)
  + \kappa \langle v \rangle ^2(\mu^{-1} \partial^\alpha_\beta f)
= -\sum_j \beta^j (\mu^{-1} \partial^{\alpha+e_j}_{\beta-e_j} f)
+ \mu^{-1} \partial^\alpha_\beta Q(f,f).
\end{align}

Set $F=\mu^{-1}f$ and denote
$F^{(\alpha,\beta)}=\mu^{-1} \partial^\alpha_\beta f=\mu^{-1} \partial_x^\alpha \partial_v^\beta f$
for $\alpha,\beta \in \mathbb{Z}^3_+$. We can rewrite Eq.(\ref{eq1}) as follows:
\begin{align}
    &\partial_t F^{(\alpha,\beta)} + v\cdot \nabla_x F^{(\alpha,\beta)} + \kappa \langle v \rangle^2 F^{(\alpha,\beta)} \\
    =&-\sum_j \beta^j (\mu^{-1} \partial^{\alpha+e_j}_{\beta-e_j} f) + \mu^{-1} \partial^\alpha_\beta Q(f,f)  \nonumber \\
    \triangleq &-\beta \cdot F^{(\alpha+1,\beta-1)}
    + \sum_{\overset{\alpha=\alpha_1+\alpha_2}{\beta=\beta_1+\beta_2}}C^{\alpha_1}_\alpha C^{\beta_1}_\beta \mu^{-1} Q(\partial^{\alpha_1}_{\beta_1} f,\partial^{\alpha_2}_{\beta_2} f)   \nonumber\\
    =&-\beta \cdot F^{(\alpha+1,\beta-1)}
     + \mu^{-1} Q(f,\partial^{\alpha}_{\beta} f)
     + \sum_{|\alpha_1|+|\beta_1| \geq 1}C^{\alpha_1}_\alpha C^{\beta_1}_\beta \mu^{-1} Q(\partial^{\alpha_1}_{\beta_1} f,\partial^{\alpha_2}_{\beta_2} f).\nonumber
\end{align}

Noticing that $\mu \mu_* = \mu' \mu'_*$, we then get the following formula,
$$
\mu^{-1}Q(f,g)=Q(\mu F,G) + \iint\! B(\mu_*-\mu'_*)F'_* G'dv_* d\sigma.
$$

Thus we obtain,
\begin{align}
   & \partial_t F^{(\alpha,\beta)} + v\cdot \nabla_x F^{(\alpha,\beta)}
      + \kappa \langle v \rangle^2 F^{(\alpha,\beta)}  \\
 = & - \beta \cdot F^{(\alpha+1,\beta-1)} + Q(\mu F,F^{(\alpha,\beta)})
     + \sum_{|\alpha_1|+|\beta_1| \geq 1} C^{\alpha_1}_\alpha C^{\beta_1}_\beta
       Q(\mu F^{(\alpha_1,\beta_1)},F^{(\alpha_2,\beta_2)})  \nonumber\\
   & + \sum_{\alpha,\beta} C^{\alpha_1}_\alpha C^{\beta_1}_\beta
     \iint B (\mu_*-\mu'_*) (F^{(\alpha_1,\beta_1)})'_* (F^{(\alpha_2,\beta_2)})'dv_* d\sigma.\nonumber
\end{align}

Set $W_l=\langle v \rangle^l $. Taking the $L^2_{x,v}$ inner product with $W_l^2 F^{(\alpha,\beta)}$, we have,
\begin{align}\label{EQ}
&\frac{1}{2}\frac{d}{dt}\|W_l F^{(\alpha,\beta)}\|^2 + \kappa \|W_{l+1} F^{(\alpha,\beta)}\|^2
 + \iint v \cdot \nabla_x F^{(\alpha,\beta)} ~ W_l^2 F^{(\alpha,\beta)} dxdv\\
=&\iint -\beta \cdot F^{(\alpha+1,\beta-1)} ~ W_l^2 F^{(\alpha,\beta)} dxdv
   + \Big( Q(\mu F,F^{(\alpha,\beta)}), W_l^2 F^{(\alpha,\beta)}\Big) \nonumber \\
 & +\sum_{|\alpha_1|+|\beta_1| \geq 1}C^{\alpha_1}_\alpha C^{\beta_1}_\beta
    \Big(Q\big(\mu F^{(\alpha_1,\beta_1)},F^{(\alpha_2,\beta_2)}\big),W_l^2 F^{(\alpha,\beta)}\Big)\nonumber\\
&+ \sum_{\alpha,\beta} C^{\alpha_1}_\alpha C^{\beta_1}_\beta
     \iiiint B (\mu_*-\mu'_*) (F^{(\alpha_1,\beta_1)})'_* (F^{(\alpha_2,\beta_2)})' W_l^2 F^{(\alpha,\beta)} dxdv dv_* d\sigma \nonumber \\
\triangleq
 &\mathcal{R}(t) + \int \Psi_1^{(0,0)}(t,x)dx
  +\sum_{|\alpha_1|+|\beta_1| \geq 1}C^{\alpha_1}_\alpha C^{\beta_1}_\beta \int \Psi_1^{(\alpha_1,\beta_1)}(t,x)dx
 +\sum_{\alpha,\beta}C^{\alpha_1}_\alpha C^{\beta_1}_\beta \int \Psi_2^{(\alpha_1,\beta_1)}(t,x)dx
  \nonumber\\
\triangleq
 &\mathcal{R}(t) + \Psi(t) + \mathcal{J}(t) + \mathcal{K}(t). \nonumber
\end{align}

We notice that the third term on the left-hand side equals to $0$,
i.e.
\begin{align}
  \iint_{{\mathbb{T}_x^3} \times {\mathbb{R}_v^3}} v \cdot \nabla_x F^{(\alpha,\beta)} ~ W_l^2
    F^{(\alpha,\beta)} dxdv
= \int_{\mathbb{R}_v^3}  v \cdot \left( \int_{\mathbb{T}_x^3} \nabla_x [\big( W_l F^{(\alpha,\beta)} \big)^2] dx  \right) dv
= 0 .
\end{align}

Multiplying by $\frac{\rho^{2(|\alpha|+|\beta|)}}{\{(\alpha-r)!\}^{2\nu_1}\{(\beta-r)!\}^{2\nu_2}}$ both sides, and integrating in time from $0$ to $t\in (0,T]$, we obtain finally,
\begin{align}\label{EQ2}
&\|f(t)\|^2_{\delta -\kappa t,l,\rho,r,\alpha,\beta}
+ 2\kappa \!\int_0^t \! \|f(\tau)\|^2_{\delta-\kappa \tau,l+1,\rho,r,\alpha,\beta}d\tau \\
\leq
&\|f(0)\|^2_{\delta,l,\rho,r,\alpha,\beta}
 +2\int^t_0\! \frac{\rho^{2(|\alpha|+|\beta|)}}{\{(\alpha-r)!\}^{2\nu_1}\{(\beta-r)!\}^{2\nu_2}}
   \Big( \mathcal{R}(\tau) + \Psi(\tau) + \mathcal{J}(\tau) + \mathcal{K}(\tau) \Big) d\tau.
   \nonumber
\end{align}

\vspace*{5pt}
\subsection{Estimates on $\mathcal{R}(t),~\Psi(t),~\mathcal{J}(t),~\mathcal{K}(t)$.}

\noindent{\bf Step I:}

Firstly we have
\begin{align}
  |\mathcal{R}(t)|
=& \left| \iint -\beta \cdot F^{(\alpha+1,\beta-1)} W_l^2 F^{(\alpha,\beta)} dxdv \right| \\
\leq & |\beta| \|W_l F^{(\alpha+1,\beta-1)}\|_{L^2_{x,v}} \!\!\!\cdot \|W_l F^{(\alpha,\beta)}\|_{L^2_{x,v}}. \nonumber
\end{align}
This yields
\begin{align}\label{R(t)}
 & \int_0^t \frac{\rho^{2(|\alpha|+|\beta|)}}{\{(\alpha-r)!\}^{2\nu_1}\{(\beta-r)!\}^{2\nu_2}}
             |\mathcal{R}(\tau)| d\tau \\
\leq & |\beta| \frac{(\alpha+1-r)^{\nu_1}}{(\beta-r)^{\nu_2}} \cdot
  \int_0^t \frac{\rho^{(|\alpha|+1)+(|\beta|-1)} \|W_l F^{(\alpha+1,\beta-1)}\|_{L^2} }
              {\{(\alpha+1-r)!\}^{\nu_1}\{(\beta-1-r)!\}^{\nu_2}}  \cdot
           \frac{\rho^{|\alpha|+|\beta|} \|W_l F^{(\alpha,\beta)}\|_{L^2} }
              {\{(\alpha-r)!\}^{\nu_1}\{(\beta-r)!\}^{\nu_2}} d\tau \nonumber \\
\leq &C_{_N} \int_0^t  \| f \|_{\delta-\kappa \tau,l,\rho,r,\alpha+1,\beta-1}
                       \| f \|_{\delta-\kappa \tau,l,\rho,r,\alpha,\beta} d\tau \nonumber \\
\leq &C \int_0^t
    (\sup_{8r\le |\alpha|+|\beta| \le N} \| f \|_{\delta-\kappa \tau,l,\rho,r,\alpha+1,\beta-1})\
    \| f \|_{l,\rho,r,N} d\tau \nonumber \\\nonumber
\le &C \int_0^t
    (\| f \|_{l,\rho,r,N} + \sup_{8r\le |\alpha| \le N}\| f \|_{\delta-\kappa \tau,l,\rho,r,\alpha+1,0})\
    \| f \|_{l,\rho,r,N} d\tau \nonumber \\\nonumber
\le &C \int_0^t (\| f \|^2_{l,\rho,r,N}
     + \sup_{8r\le |\alpha| \le N}\| f \|^2_{\delta-\kappa \tau,l,\rho,r,\alpha+1,0})d\tau, \nonumber
\end{align}
where we have used the Cauchy-Schwartz inequality in the last inequality.

\noindent{\bf Step II:}

As for $\Psi_2^{(\alpha_1,\beta_1)}(t,x)$, we recall the result in our preceding paper \cite{Zhang-Yin}:

Under the assumption $-1<\gamma+2s<1$, we have
\begin{align}
       |\Psi_2^{(\alpha_1,\beta_1)}(t,x)|
\leq C \|W_l F^{(\alpha_1,\beta_1)}\|_{L^2_v} \|W_l F^{(\alpha_2,\beta_2)}\|_{L^2_v}
         \|W_{l+1} F^{(\alpha,\beta)}\|_{L^2_v}, \,\, if \,\, l\geq 4.
\end{align}

Here we mention that, the corresponding term in \cite{Zhang-Yin} is $\Psi_2^{(\alpha',\alpha'')}(t)$, and we can consider the space variable $x$ as parameter.

Then if $|\alpha_1|+|\beta_1| \le \left[ \frac{N}{2} \right]$, we have $|\alpha_1|+|\beta_1|+2 \le N$ and
\begin{align}
\left| \int_{\mathbb{T}^3} \Psi_2^{(\alpha_1,\beta_1)}(t,x) dx \right|
\leq & C \|W_l F^{(\alpha_1,\beta_1)}\|_{L_x^{\infty}(L^2_v)}
         \|W_l F^{(\alpha_2,\beta_2)}\|_{L_x^2(L^2_v)}
         \|W_{l+1} F^{(\alpha,\beta)}\|_{L_x^2(L^2_v)} \\
\leq & C \|W_l F^{(\alpha_1,\beta_1)}\|_{H_x^2(L^2_v)}
         \|W_l F^{(\alpha_2,\beta_2)}\|_{L_x^2(L^2_v)}
         \|W_{l+1} F^{(\alpha,\beta)}\|_{L_x^2(L^2_v)} \nonumber\\
\leq & C \|W_l F^{(\alpha_1+2,\beta_1)}\|_{L^2_{x,v}} \|W_l F^{(\alpha_2,\beta_2)}\|_{L^2_{x,v}}
         \|W_{l+1} F^{(\alpha,\beta)}\|_{L^2_{x,v}}, \nonumber
\end{align}
where in the second inequality we have used the embedding $\|h\|_{H_x^2} \leq C \|h\|_{L_x^\infty}$.

On the other hand, if $|\alpha_1|+|\beta_1| \ge \left[ \frac{N}{2} \right]+1$, which implies that $|\alpha_2|+|\beta_2|+2 \le N$, we have
\begin{align}
\left| \int_{\mathbb{T}^3} \Psi_2^{(\alpha_1,\beta_1)}(t,x) dx \right|
\leq & C \|W_l F^{(\alpha_1,\beta_1)}\|_{L_x^2(L^2_v)}
         \|W_l F^{(\alpha_2,\beta_2)}\|_{L_x^{\infty}(L^2_v)}
         \|W_{l+1} F^{(\alpha,\beta)}\|_{L_x^2(L^2_v)} \\
\leq & C \|W_l F^{(\alpha_1,\beta_1)}\|_{L_x^2(L^2_v)}
         \|W_l F^{(\alpha_2,\beta_2)}\|_{H_x^2(L^2_v)}
         \|W_{l+1} F^{(\alpha,\beta)}\|_{L_x^2(L^2_v)} \nonumber\\
\leq & C \|W_l F^{(\alpha_1,\beta_1)}\|_{L^2_{x,v}} \|W_l F^{(\alpha_2+2,\beta_2)}\|_{L^2_{x,v}}
         \|W_{l+1} F^{(\alpha,\beta)}\|_{L^2_{x,v}}. \nonumber
\end{align}

Therefore, we obtain if $|\alpha_1|+|\beta_1| \le \left[ \frac{N}{2} \right]$, then
\begin{align}\label{Psi_2}
\quad \frac{\rho^{2(|\alpha|+|\beta|)}
                   \left| \int \Psi_2^{(\alpha_1,\beta_1)}(t,x) dx \right|}
                 {\{(\alpha-r)!\}^{2\nu_1}\{(\beta-r)!\}^{2\nu_2}}
\leq & C \frac{\{(\alpha_1+2-r)!\}^{\nu_1}
\{(\alpha_2-r)!\}^{\nu_1}} {\{(\alpha-r)!\}^{\nu_1}} \cdot
         \frac{\{(\beta_1-r)!\}^{\nu_2} \{(\beta_2-r)!\}^{\nu_2} } {\{(\beta-r)!\}^{\nu_2}} \\
       &   \times \| f \|_{\delta-\kappa t,l,\rho,r,\alpha_1+2,\beta_1} ~
                  \| f \|_{\delta-\kappa t,l,\rho,r,\alpha_2,\beta_2} ~
                  \| f \|_{\delta-\kappa t,l+1,\rho,r,\alpha,\beta} ~ ,
\nonumber
\end{align}
and if $|\alpha_1|+|\beta_1| \ge \left[ \frac{N}{2} \right]+1$, then
\begin{align}\label{Psi'_2}
\quad \frac{\rho^{2(|\alpha|+|\beta|)}
                   \left| \int \Psi_2^{(\alpha_1,\beta_1)}(t,x) dx \right|}
                 {\{(\alpha-r)!\}^{2\nu_1}\{(\beta-r)!\}^{2\nu_2}}
\leq & C \frac{\{(\alpha_1-r)!\}^{\nu_1}
\{(\alpha_2+2-r)!\}^{\nu_1}} {\{(\alpha-r)!\}^{\nu_1}} \cdot
         \frac{\{(\beta_1-r)!\}^{\nu_2} \{(\beta_2-r)!\}^{\nu_2} } {\{(\beta-r)!\}^{\nu_2}} \\
       &   \times \| f \|_{\delta-\kappa t,l,\rho,r,\alpha_1,\beta_1} ~
                  \| f \|_{\delta-\kappa t,l,\rho,r,\alpha_2+2,\beta_2} ~
                  \| f \|_{\delta-\kappa t,l+1,\rho,r,\alpha,\beta} ~ .
\nonumber
\end{align}

Specifically, we point out that if $\beta=0$, i.e., $(\alpha,\ \beta)=(\alpha,\ 0)$ and $\beta_1=\beta_2=0$, the above two estimates remain valid.

\noindent{\bf Step III:}

Next, concerning $\Psi_1^{(0,0)}(t,x)$, we use the coercivity estimate in \cite{Zhang-Yin}:
\begin{align}
\Psi_1^{(0,0)}(t,x)+ c_0 \| W_{l+\gamma/2} F^{(\alpha,\beta)} \|^2_{H^s_v}
&\leq C \|W_{l+\gamma} F^{(\alpha,\beta)}\|_{L^2_v} \| W_l F^{(\alpha,\beta)} \|_{L^2_v} \\
&\leq C \|W_{l+1} F^{(\alpha,\beta)}\|_{L^2_v} \| W_l
F^{(\alpha,\beta)} \|_{L^2_v} \,. \nonumber
\end{align}
Then it follows
\begin{align}
   \int_{\mathbb{T}^3} \Psi_1^{(0,0)}(t,x)dx
    + c_0 \| W_{l+\gamma/2} F^{(\alpha,\beta)} \|^2_{L^2_x (H^s_v)}
  \leq C \|W_l F^{(\alpha,\beta)}\|_{L^2_{x,v}} \| W_{l+1} F^{(\alpha,\beta)} \|_{L^2_{x,v}}.
\end{align}

Observing that $\Psi(t)=\int_{\mathbb{T}^3} \Psi_1^{(0,0)}(t,x)dx$, we have
\begin{align}\label{Psi(t)}
 & \int_0^t \frac{\rho^{2(|\alpha|+|\beta|)} \,\,\Psi(\tau) }
                 {\{(\alpha-r)!\}^{2\nu_1}\{(\beta-r)!\}^{2\nu_2}} d\tau
+c_0 \int_0^t \frac{\rho^{2(|\alpha|+|\beta|)}
                     \| W_{l+\gamma/2} F^{(\alpha,\beta)} \|^2_{L^2_x (H^s_v)}}
                   {\{(\alpha-r)!\}^{2\nu_1}\{(\beta-r)!\}^{2\nu_2}} d\tau \\
\leq & C \, \int_0^t \| f \|_{\delta-\kappa \tau,l,\rho,r,\alpha,\beta} \,
                     \| f \|_{\delta-\kappa \tau,l+1,\rho,r,\alpha,\beta} d\tau \nonumber \\
\leq & C_\kappa \, \int_0^t \| f \|^2_{l,\rho,r,N} d\tau \,
         +\frac{\kappa}{100} \sup_{8r \leq |\alpha|+|\beta| \leq N}
            \int_0^t \| f \|^2_{\delta-\kappa \tau,l+1,\rho,r,\alpha,\beta} d\tau. \nonumber
\end{align}

\noindent{\bf Step IV:}

Now we consider the term $\Psi_1^{(\alpha_1,\beta_1)}(t,x)$ with the restriction $|\alpha_1|+|\beta_1| \geq 1$. By virtue of the upper bound estimate for collision operator $Q$ from Proposition 2.9 of \cite{FiveGroup-III}:
\begin{lemm}\label{upper bound-2}
Let $0<s<1$ and $-1<\gamma+2s < 1$. For any $p\in \mathbb{R}$ and $m\in [s-1,s]$, there exists a $C>0$ such that
\begin{align*}
\Big| \Big( Q(f,g),h \Big)_{L^2(\mathbb{R}^3)} \Big |
\leq C \big( \|f\|_{L^1_{p^++(\gamma+2s)^+}} + \|f\|_{L^2} \big)
\|g\|_{H^{\max{\{s+m,(2s-1+\varepsilon)^+\}}}_{(p+\gamma+2s)^+}} \|h\|_{H^{s-m}_{-p}}\,.
\end{align*}
\end{lemm}

Applying above lemma with $f=\mu F^{(\alpha_1,\beta_1)}$, $g=F^{(\alpha_2,\beta_2)}$, $h=W_{2l} F^{(\alpha,\beta)}$, $p=l-\gamma-2s$, and $m=s$, by setting $\eta=1-(\gamma+2s) \text{ if }\gamma+2s \in (0,1)$, we can infer that,
\begin{align*}
\left| \Psi_1^{(\alpha_1,\beta_1)} \right|
\leq &C
\big(\| \mu F^{(\alpha_1,\beta_1)} \|_{L^1_{l-\gamma -2s+(\gamma+2s)^+}}
     + \| \mu F^{(\alpha_1,\beta_1)} \|_{L^2} \big)
\| W_l F^{(\alpha_2,\beta_2)} \|_{H^{2s}}
\| W_{l+\gamma +2s} F^{(\alpha,\beta)} \|_{L^2}\\
\leq &C\| W_l F^{(\alpha_1,\beta_1)} \|_{L^2} \|W_{l+1} F^{(\alpha_2,\beta_2)} \|_{L^2} \| W_{l+1-\eta} F^{(\alpha,\beta)} \|_{L^2} \\
& + C \| W_l F^{(\alpha_1,\beta_1)} \|_{L^2} \|W_l F^{(\alpha_2,\beta_2+1)} \|_{L^2} \| W_{l+1} F^{(\alpha,\beta)} \|_{L^2},
\end{align*}
in view of $2s<1$ and
$$
\partial_v (W_l \mu^{-1} \partial^{\alpha_2}_{\beta_2}f)=\partial_v (W_l \mu^{-1}) \partial^{\alpha_2}_{\beta_2}f
+ W_l\mu^{-1} \partial^{\alpha_2}_{\beta_2+1}f.
$$

Whereas in the case $\gamma+2s \in (-1,0]$, by setting $\eta=1+(\gamma+2s) \in (0,1]$, we have
\begin{align*}
\left| \Psi_1^{(\alpha_1,\beta_1)} \right|
\leq &C
\big(\| \mu F^{(\alpha_1,\beta_1)} \|_{L^1_{l-\gamma -2s+(\gamma+2s)^+}}
     + \| \mu F^{(\alpha_1,\beta_1)} \|_{L^2} \big)
\| W_l F^{(\alpha_2,\beta_2)} \|_{H^{2s}}
\| W_{l+\gamma +2s} F^{(\alpha,\beta)} \|_{L^2}\\
\leq &C\| W_l F^{(\alpha_1,\beta_1)} \|_{L^2} \|W_{l+1} F^{(\alpha_2,\beta_2)} \|_{L^2} \| W_{l} F^{(\alpha,\beta)} \|_{L^2} \\
& + C \| W_l F^{(\alpha_1,\beta_1)} \|_{L^2} \|W_l F^{(\alpha_2,\beta_2+1)} \|_{L^2} \| W_{l} F^{(\alpha,\beta)} \|_{L^2} \\
\leq &C\| W_l F^{(\alpha_1,\beta_1)} \|_{L^2} \|W_{l+1} F^{(\alpha_2,\beta_2)} \|_{L^2} \| W_{l+1-\eta} F^{(\alpha,\beta)} \|_{L^2} \\
& + C \| W_l F^{(\alpha_1,\beta_1)} \|_{L^2} \|W_l F^{(\alpha_2,\beta_2+1)} \|_{L^2} \| W_{l+1} F^{(\alpha,\beta)} \|_{L^2}.
\end{align*}

Note that we have the same estimate on $\Psi_1^{(\alpha_1,\beta_1)}$ in both cases $0<\gamma+2s<1$ and $-1<\gamma+2s \leq 0$, then if $|\alpha_1|+|\beta_1| \le \left[ \frac{N}{2} \right]$, by recalling $|\alpha_1|+|\beta_1|\ge 1$, it holds $|\alpha_2|+|\beta_2|+1 \le N$, hence we can write that
\begin{align}
\left| \int \Psi_1^{(\alpha_1,\beta_1)}(t,x) dx \right|
\leq &C\|W_l F^{(\alpha_1,\beta_1)}\|_{L_x^\infty(L^2_v)}
       \|W_{l+1} F^{(\alpha_2,\beta_2)}\|_{L_x^2(L^2_v)}
       \|W_{l+1-\eta} F^{(\alpha,\beta)}\|_{L_x^2(L^2_v)} \\
    &+ C\|W_l F^{(\alpha_1,\beta_1)}\|_{L_x^\infty(L^2_v)}
        \|W_l F^{(\alpha_2,\beta_2+1)}\|_{L_x^2(L^2_v)}
        \|W_{l+1} F^{(\alpha,\beta)}\|_{L_x^2(L^2_v)}  \nonumber \\
\leq &C\|W_l F^{(\alpha_1+2,\beta_1)}\|_{L^2_{x,v}}
       \|W_{l+1} F^{(\alpha_2,\beta_2)}\|_{L^2_{x,v}}
       \|W_l F^{(\alpha,\beta)}\|_{L^2_{x,v}}\nonumber \\
    &+ C\|W_l F^{(\alpha_1+2,\beta_1)}\|_{L^2_{x,v}}
        \|W_l F^{(\alpha_2,\beta_2+1)}\|_{L^2_{x,v}}
        \|W_{l+1} F^{(\alpha,\beta)}\|_{L^2_{x,v}}   \nonumber
\end{align}

Since the H\"{o}lder inequality yields
$$
\| W_{1-\eta} G\|_{L^2}\leq \|G\|^\eta_{L^2} \|W_1 G\|^{1-\eta}_{L^2},
$$
we obtain if $|\alpha_1|+|\beta_1| \le \left[ \frac{N}{2} \right]$,
\begin{align}
\left| \int_{\mathbb{T}^3} \Psi_1^{(\alpha_1,\beta_1)}(t,x) dx \right|
\leq & C \|W_l F^{(\alpha_1+2,\beta_1)}\|_{L^2_{x,v}} \|W_{l+1} F^{(\alpha_2,\beta_2)}\|_{L^2_{x,v}}
    \|W_{l+1} F^{(\alpha,\beta)}\|^{1-\eta}_{L^2_{x,v}} \|W_l F^{(\alpha,\beta)}\|^\eta_{L^2_{x,v}} \\
     & + C\|W_l F^{(\alpha_1+2,\beta_1)}\|_{L^2_{x,v}} \|W_l F^{(\alpha_2,\beta_2+1)}\|_{L^2_{x,v}}
          \|W_{l+1} F^{(\alpha,\beta)}\|_{L^2_{x,v}} \nonumber \\
\triangleq & J_1(t)+J_2(t). \nonumber
\end{align}

Or alternatively, when $|\alpha_1|+|\beta_1| \ge \left[ \frac{N}{2} \right]+1$, we can arrive at $|\alpha_2|+|\beta_2|+2+1 \le N-(\left[ \frac{N}{2} \right]+1) +2+1 \le N$, and further,
\begin{align}
\left| \int_{\mathbb{T}^3} \Psi_1^{(\alpha_1,\beta_1)}(t,x) dx \right|
\leq & C \|W_l F^{(\alpha_1,\beta_1)}\|_{L^2_{x,v}} \|W_{l+1} F^{(\alpha_2+2,\beta_2)}\|_{L^2_{x,v}}
    \|W_{l+1} F^{(\alpha,\beta)}\|^{1-\eta}_{L^2_{x,v}} \|W_l F^{(\alpha,\beta)}\|^\eta_{L^2_{x,v}} \\
     & + C\|W_l F^{(\alpha_1,\beta_1)}\|_{L^2_{x,v}} \|W_l F^{(\alpha_2+2,\beta_2+1)}\|_{L^2_{x,v}}
          \|W_{l+1} F^{(\alpha,\beta)}\|_{L^2_{x,v}} \nonumber \\
\triangleq & J_1(t)+J_2(t). \nonumber
\end{align}
Above, $\eta=1-(\gamma+2s) \text{ if }\gamma+2s \in (0,1)$, or $\eta=1+(\gamma+2s)$ if $\gamma+2s \in (-1,0]$.

Thus, if $|\alpha_1|+|\beta_1| \le \left[ \frac{N}{2} \right]$, we get that
\begin{align}\label{J_1-original}
  \frac{\rho^{2(|\alpha|+|\beta|)} J_1(t)}{\{(\alpha-r)!\}^{2\nu_1}\{(\beta-r)!\}^{2\nu_2}}
\leq & C \frac{\{(\alpha_1+2-r)!\}^{\nu_1}
\{(\alpha_2-r)!\}^{\nu_1}} {\{(\alpha-r)!\}^{\nu_1}} \cdot
         \frac{\{(\beta_1-r)!\}^{\nu_2} \{(\beta_2-r)!\}^{\nu_2} } {\{(\beta-r)!\}^{\nu_2}} \\
     &  \times \| f \|_{\delta-\kappa t,l,\rho,r,\alpha_1+2,\beta_1} ~
               \| f \|^\eta_{\delta-\kappa t,l,\rho,r,\alpha,\beta} ~
               \| f \|_{\delta-\kappa t,l+1,\rho,r,\alpha_2,\beta_2} ~
               \| f \|^{1-\eta}_{\delta-\kappa t,l+1,\rho,r,\alpha,\beta}~, \nonumber
\end{align}
and
\begin{align}\label{J_2-original}
\frac{\rho^{2(|\alpha|+|\beta|)} J_2(t)}{\{(\alpha-r)!\}^{2\nu_1}\{(\beta-r)!\}^{2\nu_2}}
\leq & C \frac{\{(\alpha_1+2-r)!\}^{\nu_1}
\{(\alpha_2-r)!\}^{\nu_1}} {\{(\alpha-r)!\}^{\nu_1}} \cdot
         \frac{\{(\beta_1-r)!\}^{\nu_2} \{(\beta_2+1-r)!\}^{\nu_2} } {\{(\beta-r)!\}^{\nu_2}} \\
     &  \times \| f \|_{\delta-\kappa t,l,\rho,r,\alpha_1+2,\beta_1} ~
               \| f \|_{\delta-\kappa t,l,\rho,r,\alpha_2,\beta_2+1} ~
               \| f \|_{\delta-\kappa t,l+1,\rho,r,\alpha,\beta} ~.\nonumber
\end{align}

Similarly, when $|\alpha_1|+|\beta_1| \ge \left[ \frac{N}{2} \right]+1$, we can have
\begin{align}\label{J'_1-original}
  \frac{\rho^{2(|\alpha|+|\beta|)} J_1(t)}{\{(\alpha-r)!\}^{2\nu_1}\{(\beta-r)!\}^{2\nu_2}}
\leq & C \frac{\{(\alpha_1-r)!\}^{\nu_1}
\{(\alpha_2+2-r)!\}^{\nu_1}} {\{(\alpha-r)!\}^{\nu_1}} \cdot
         \frac{\{(\beta_1-r)!\}^{\nu_2} \{(\beta_2-r)!\}^{\nu_2} } {\{(\beta-r)!\}^{\nu_2}} \\
     &  \times \| f \|_{\delta-\kappa t,l,\rho,r,\alpha_1,\beta_1} ~
               \| f \|^\eta_{\delta-\kappa t,l,\rho,r,\alpha,\beta} ~
               \| f \|_{\delta-\kappa t,l+1,\rho,r,\alpha_2+2,\beta_2} ~
               \| f \|^{1-\eta}_{\delta-\kappa t,l+1,\rho,r,\alpha,\beta}~, \nonumber
\end{align}
and
\begin{align}\label{J'_2-original}
\frac{\rho^{2(|\alpha|+|\beta|)} J_2(t)}{\{(\alpha-r)!\}^{2\nu_1}\{(\beta-r)!\}^{2\nu_2}}
\leq & C \frac{\{(\alpha_1-r)!\}^{\nu_1}
\{(\alpha_2+2-r)!\}^{\nu_1}} {\{(\alpha-r)!\}^{\nu_1}} \cdot
         \frac{\{(\beta_1-r)!\}^{\nu_2} \{(\beta_2+1-r)!\}^{\nu_2} } {\{(\beta-r)!\}^{\nu_2}} \\
     &  \times \| f \|_{\delta-\kappa t,l,\rho,r,\alpha_1,\beta_1} ~
               \| f \|_{\delta-\kappa t,l,\rho,r,\alpha_2+2,\beta_2+1} ~
               \| f \|_{\delta-\kappa t,l+1,\rho,r,\alpha,\beta} ~.\nonumber
\end{align}

\noindent{\bf Step V:}

Before continuing the proof, we state a useful lemma, as follows:
\begin{lemm}
If $\nu \geq 1$ and $2\leq r \in \mathbb{N}$, then there exists a constant $B>0$ depending only on $r$ such that for any $\alpha \in \mathbb{Z}^3$,
\begin{align}
  \sum_{\alpha=\alpha_1+\alpha_2} C_\alpha^{\alpha_1}
      \frac{\{(\alpha_1-r)!\}^\nu\{(\alpha_2-r)!\}^\nu}{\{(\alpha-r)!\}^\nu} \leq B,
\end{align}
and
\begin{align}
  \sum_{\alpha=\alpha_1+\alpha_2} C_\alpha^{\alpha_1}
      \frac{\{(\alpha_1+2-r)!\}^{\nu}\{(\alpha_2-r)!\}^{\nu}} {\{(\alpha-r)!\}^{\nu}} \leq C \cdot N^{2\nu} \triangleq B'.
\end{align}
Furthermore, if $\nu>1$ and $r>1+\nu/(\nu-1)$, then there exists a constant $B''>0$ depending only on $\nu$ and $r$ such that for any $0 \neq \alpha \in \mathbb{Z}^3$,
\begin{align}
  \sum_{\alpha=\alpha_1+\alpha_2,\alpha_1 \neq 0} C_\alpha^{\alpha_1}
      \frac{\{(\alpha_1-r)!\}^\nu\{(\alpha_2+1-r)!\}^\nu}{\{(\alpha-r)!\}^\nu} \leq B''.
\end{align}
\end{lemm}

Note that the first and third result were established in \cite{Mori-Ukai}, and the proof of the second share a similar scheme as the first, so we omit it here. Additionally, we will consider $B$, $B'$ and $B''$ as $1$ in later proof.

Now we resume the proof of Lemma \ref{Key lemma}. As for the integral including $\mathcal{K}(t)$ in (\ref{EQ2}), from (\ref{Psi_2}) we can infer that, if $|\alpha_1|+|\beta_1| \le \left[ \frac{N}{2} \right]$,
\begin{align}
  &\int_0^t \frac{\rho^{2(|\alpha|+|\beta|)} ~ |\mathcal{K}(\tau)| }
                 {\{(\alpha-r)!\}^{2\nu_1}\{(\beta-r)!\}^{2\nu_2}}  d\tau \\\nonumber
\leq & C \sum_{\alpha,\beta} C^{\alpha_1}_\alpha C^{\beta_1}_\beta
               \frac{\{(\alpha_1+2-r)!\}^{\nu_1} \{(\alpha_2-r)!\}^{\nu_1}}
                    {\{(\alpha-r)!\}^{\nu_1}}
         \frac{\{(\beta_1-r)!\}^{\nu_2} \{(\beta_2-r)!\}^{\nu_2} } {\{(\beta-r)!\}^{\nu_2}} \\
     &  \times \int_0^t \| f \|_{\delta-\kappa \tau,l,\rho,r,\alpha_1+2,\beta_1} ~
                        \| f \|_{\delta-\kappa \tau,l,\rho,r,\alpha_2,\beta_2} ~
                        \| f \|_{\delta-\kappa \tau,l+1,\rho,r,\alpha,\beta} ~ d\tau \nonumber \\
\leq & C \int_0^t \| f \|_{\delta-\kappa \tau,l,\rho,r,\alpha_1+2,\beta_1} ~
                  \| f \|_{\delta-\kappa \tau,l,\rho,r,\alpha_2,\beta_2} ~
                  \| f \|_{\delta-\kappa \tau,l+1,\rho,r,\alpha,\beta} ~ d\tau \nonumber \\
\leq & C \frac{1}{4 \varepsilon} \int_0^t \| f \|^2_{\delta-\kappa \tau,l,\rho,r,\alpha_1+2,\beta_1}~
                        \| f \|^2_{\delta-\kappa \tau,l,\rho,r,\alpha_2,\beta_2}~ d\tau
      + \varepsilon \int_0^t \|f(\tau)\|^2_{\delta-\kappa \tau,l+1,\rho,r,\alpha,\beta}~ d\tau,\nonumber
\end{align}
and if $|\alpha_1|+|\beta_1| \ge \left[ \frac{N}{2} \right]+1$, we have, correspondingly,
\begin{align}
 & \int_0^t \frac{\rho^{2(|\alpha|+|\beta|)} ~ |\mathcal{K}(\tau)| }
                 {\{(\alpha-r)!\}^{2\nu_1}\{(\beta-r)!\}^{2\nu_2}}  d\tau \\\nonumber
\leq & C \sum_{\alpha,\beta} C^{\alpha_1}_\alpha C^{\beta_1}_\beta
               \frac{\{(\alpha_1-r)!\}^{\nu_1} \{(\alpha_2+2-r)!\}^{\nu_1}}
                    {\{(\alpha-r)!\}^{\nu_1}}
         \frac{\{(\beta_1-r)!\}^{\nu_2} \{(\beta_2-r)!\}^{\nu_2} } {\{(\beta-r)!\}^{\nu_2}} \\
     &  \times \int_0^t \| f \|_{\delta-\kappa \tau,l,\rho,r,\alpha_1,\beta_1} ~
                        \| f \|_{\delta-\kappa \tau,l,\rho,r,\alpha_2+2,\beta_2} ~
                        \| f \|_{\delta-\kappa \tau,l+1,\rho,r,\alpha,\beta} ~ d\tau \nonumber \\
\leq & C \int_0^t \| f \|_{\delta-\kappa \tau,l,\rho,r,\alpha_1,\beta_1} ~
                  \| f \|_{\delta-\kappa \tau,l,\rho,r,\alpha_2+2,\beta_2} ~
                  \| f \|_{\delta-\kappa \tau,l+1,\rho,r,\alpha,\beta} ~ d\tau \nonumber \\
\leq & C \frac{1}{4 \varepsilon} \int_0^t \| f \|^2_{\delta-\kappa \tau,l,\rho,r,\alpha_1,\beta_1}~
                        \| f \|^2_{\delta-\kappa \tau,l,\rho,r,\alpha_2+2,\beta_2}~ d\tau
      + \varepsilon \int_0^t \|f(\tau)\|^2_{\delta-\kappa \tau,l+1,\rho,r,\alpha,\beta}~ d\tau,\nonumber
\end{align}

We then introduce the following notations:
\begin{align}
  & A_l \xlongequal{\!\!def\!\!}
  \sup_{t\in [0,T]} \sup_{|\alpha|+|\beta| \leq 8r} \| f(t)\|_{\delta-\kappa t,l,\rho,r,\alpha,\beta}.
%  & [f]_{l,\rho,r,N}(t,\alpha_0) \triangleq \sup_{\overset{ |\alpha'| \leq 3r}{3r \leq |\beta'| \leq N}}
%                   \| f(t)\|_{\delta-\kappa t,l,\rho,r,\alpha',\beta'},
\end{align}
%where $N$ takes the same value as before. Besides, $\alpha_0$ only stands for some multi-index less than $3r$, and there has been no need to know its exact value. Likewise, we can also define $[f]_{l,\rho,r,N}(t,\beta_0)$.

Now, we shall give the different bounds on the factor
\begin{align}
\Xi ~ \triangleq ~  & \sup_{8r\le |\alpha|+|\beta| \leq N}
      \| f \|^2_{\delta-\kappa \tau,l,\rho,r,\alpha_1+2,\beta_1} ~
      \| f \|^2_{\delta-\kappa \tau,l,\rho,r,\alpha_2,\beta_2},\nonumber
\end{align}
with respect to the values of $|\alpha_1|+|\beta_1|$ and $|\alpha_2|+|\beta_2|$:
\begin{itemize}
  \item[$\cdot$] if $|\alpha_1|+|\beta_1|\le 8r-2,\ |\alpha_2|+|\beta_2| \leq 8r$, we have
          $$\Xi \leq \sup_{|\alpha_1|+|\beta_1|+2 \leq 8r}
                    \| f \|^2_{\delta-\kappa \tau,l,\rho,r,\alpha_1+2,\beta_1}
                    \times \sup_{|\alpha_2|+|\beta_2| \leq 8r}
                    \| f \|^2_{\delta-\kappa \tau,l,\rho,r,\alpha_2,\beta_2}
                \le A_l^4;$$
  \item[$\cdot$] if $|\alpha_1|+|\beta_1|\le 8r-2,\ |\alpha_2|+|\beta_2| \geq 8r$, we have
          $$\Xi \leq \sup_{|\alpha_1|+|\beta_1|+2 \leq 8r}
                    \| f \|^2_{\delta-\kappa \tau,l,\rho,r,\alpha_1+2,\beta_1}
                    \times \sup_{8r\le |\alpha_2|+|\beta_2| \leq N}
                    \| f \|^2_{\delta-\kappa \tau,l,\rho,r,\alpha_2,\beta_2}
                \le A_l^2 \|f\|^2_{l,\rho,r,N};$$
  \item[$\cdot$] if $|\alpha_1|+|\beta_1|\ge 8r-2,\ |\alpha_2|+|\beta_2| \leq 8r$, we have
          $$\Xi \leq \sup_{8r\le |\alpha_1|+|\beta_1|+2 \leq N}
                    \| f \|^2_{\delta-\kappa \tau,l,\rho,r,\alpha_1+2,\beta_1}
                    \times \sup_{|\alpha_2|+|\beta_2| \leq 8r}
                    \| f \|^2_{\delta-\kappa \tau,l,\rho,r,\alpha_2,\beta_2}
                \le A_l^2 \|f\|^2_{l,\rho,r,N};$$
  \item[$\cdot$] if $|\alpha_1|+|\beta_1|\ge 8r-2,\ |\alpha_2|+|\beta_2| \geq 8r$, we have
          $$\Xi \leq \sup_{8r\le |\alpha_1|+|\beta_1|+2 \leq N}
                    \| f \|^2_{\delta-\kappa \tau,l,\rho,r,\alpha_1+2,\beta_1}
                    \times \sup_{8r\le |\alpha_2|+|\beta_2| \leq N}
                    \| f \|^2_{\delta-\kappa \tau,l,\rho,r,\alpha_2,\beta_2}
                \le \|f\|^4_{l,\rho,r,N};$$
\end{itemize}

Thus, together with the above inequalities, we get,
\begin{align}
  \Xi ~ \leq ~ \|f\|^4_{l,\rho,r,N} + A_l^2 \|f\|^2_{l,\rho,r,N} + A_l^4,
\end{align}
which yields that,
\begin{align}\label{K(t)}
   \int_0^t \frac{\rho^{2(|\alpha|+|\beta|)} ~ |\mathcal{K}(\tau)| }
                 {\{(\alpha-r)!\}^{2\nu_1}\{(\beta-r)!\}^{2\nu_2}}  d\tau
 \leq & C_\kappa \int_0^t
           \left\{ \|f\|^4_{l,\rho,r,N}(\tau) + A_l^2~ \|f\|^2_{l,\rho,r,N}(\tau) + A_l^4
           \right\} d\tau \\ &
      + \frac{\kappa}{100} \sup_{8r \leq |\alpha|+|\beta| \leq N}
                         \int_0^t \|f(\tau)\|^2_{\delta-\kappa \tau,l+1,\rho,r,\alpha,\beta} d\tau.
 \nonumber
\end{align}

On the other hand, in the case $|\alpha_1|+|\beta_1| \ge \left[ \frac{N}{2} \right]+1$, a similar scheme ensures us to get the same estimate as above.

Considering the integral including $\mathcal{J}(t)$ in (\ref{EQ2}), we can write that,
\begin{align*}
\int_0^t \frac{\rho^{2(|\alpha|+|\beta|)} \,\,|\mathcal{J}(\tau)| }
                 {\{(\alpha-r)!\}^{2\nu_1}\{(\beta-r)!\}^{2\nu_2}}  d\tau
\leq & \int_0^t \sum_{\alpha,\beta} C^{\alpha_1}_\alpha C^{\beta_1}_\beta
                \frac{\rho^{2(|\alpha|+|\beta|)} \,\,\mathcal{J}_1(\tau) }
                 {\{(\alpha-r)!\}^{2\nu_1}\{(\beta-r)!\}^{2\nu_2}}  d\tau \\
    & + \int_0^t \sum_{|\alpha_1|+|\beta_1| \geq 1}
                        C^{\alpha_1}_\alpha C^{\beta_1}_\beta
                    \frac{\rho^{2(|\alpha|+|\beta|)} \,\,\mathcal{J}_2(\tau) }
                     {\{(\alpha-r)!\}^{2\nu_1}\{(\beta-r)!\}^{2\nu_2}}  d\tau.
\end{align*}

First, it is easy to check that, from (\ref{J_2-original}) and (\ref{J'_2-original}), the integral including $\mathcal{J}_2(t)$ has the same estimate as the inequality (\ref{K(t)}).

Furthermore, the last factor of (\ref{J_1-original}) (or (\ref{J'_1-original})) is bounded by
\begin{align}
  C_\varepsilon \|f(t)\|^{2/\eta}_{\delta-\kappa t,l,\rho,r,\alpha_1+2,\beta_1}
                        \|f(t)\|^2_{\delta-\kappa t,l,\rho,r,\alpha,\beta}
   + \varepsilon \Big( \|f(t)\|^2_{\delta-\kappa t,l+1,\rho,r,\alpha_2,\beta_2}
                        + \|f(t)\|^2_{\delta-\kappa t,l+1,\rho,r,\alpha,\beta}
                  \Big).
\end{align}
or
\begin{align}
  C_\varepsilon \|f(t)\|^{2/\eta}_{\delta-\kappa t,l,\rho,r,\alpha_1,\beta_1}
                        \|f(t)\|^2_{\delta-\kappa t,l,\rho,r,\alpha,\beta}
   + \varepsilon \Big( \|f(t)\|^2_{\delta-\kappa t,l+1,\rho,r,\alpha_2+2,\beta_2}
                        + \|f(t)\|^2_{\delta-\kappa t,l+1,\rho,r,\alpha,\beta}
                  \Big).
\end{align}

A similar procedure as that of handling with the estimate of (\ref{K(t)}) gives that,
\begin{align}
  \int_0^t \frac{\rho^{2(|\alpha|+|\beta|)} \,\,\mathcal{J}_1(\tau) }
                 {\{(\alpha-r)!\}^{2\nu_1}\{(\beta-r)!\}^{2\nu_2}}  d\tau
\leq & C_\kappa \int_0^t \left\{ \|f\|^{2(1+\eta)/\eta}_{l,\rho,r,N}(\tau)
                               + A_l^{2/\eta} ~ \|f\|^2_{l,\rho,r,N}(\tau)
                         \right\} d\tau \\
    & + \frac{\kappa}{100}
          \sup_{8r \leq |\alpha|+|\beta| \leq N}
            \int_0^t \Big( \|f(\tau)\|^2_{\delta-\kappa \tau,l+1,\rho,\alpha,\beta} + A_{l+1}^2
                     \Big) d\tau. \nonumber
\end{align}

Observing that $0<\eta<1$ implies $2/\eta >2$ and $2(1+\eta)/\eta >4$, so we obtain the estimate on $\mathcal{J}(t)$ as follows:
\begin{align}\label{J(t)}
   \int_0^t \frac{\rho^{2(|\alpha|+|\beta|)} ~ |\mathcal{J}(\tau)| }
                 {\{(\alpha-r)!\}^{2\nu_1}\{(\beta-r)!\}^{2\nu_2}}  d\tau
\leq & C_\kappa \int_0^t \left\{ \|f\|^{2(1+\eta)/\eta}_{l,\rho,r,N}(\tau)
                               + A_l^{2/\eta}\, \|f\|^2_{l,\rho,r,N}(\tau) + A_l^4
                        \right\} d\tau \\
    & + \frac{\kappa}{50} \sup_{8r \leq |\alpha|+|\beta| \leq N}
                         \int_0^t \Big( \|f(\tau)\|^2_{\delta-\kappa \tau,l+1,\rho,\alpha,\beta} + A_{l+1}^2
                                  \Big) d\tau.\nonumber
\end{align}

\subsection{Main inequalities and their proofs}

\noindent{\bf Step I:}
Combining (\ref{EQ2}), (\ref{R(t)}), (\ref{Psi(t)}), (\ref{K(t)}) and (\ref{J(t)}), we finally deduce that,
\begin{align}\label{Gronw}
  &\| f(t)\|^2_{\delta-\kappa t,l,\rho,r,\alpha,\beta}
 + 2 c_0 \int_0^t
      \frac{\rho^{2(|\alpha|+|\beta|)} \| W_{l+\frac{\gamma}{2}} \mu^{-1} \partial^\alpha_\beta f(\tau)\|^2_{L^2_x(H^s_v)}} {\{ (\alpha-r)!\}^{2\nu_1} \{ (\beta-r)!\}^{2\nu_2}} d\tau
 + 2\kappa \int_0^t \| f(\tau) \|^2_{\delta-\kappa \tau,l+1,\rho,r,\alpha,\beta} d\tau  \\\nonumber
\leq & \|f(0)\|^2_{\delta,l,\rho,r,\alpha,\beta}
  + C_\kappa \int_0^t \left\{  \|f\|^{2(1+\eta)/\eta}_{l,\rho,r,N}(\tau)
                             + A_l^{2/\eta} \|f \|^2_{l,\rho,r,N}(\tau) + A_l^4
                             + \sup_{8r\le |\alpha| \le N} \| f \|^2_{\delta-\kappa \tau,l,\rho,r,\alpha+1,0}(\tau)
                      \right\} d\tau \\\nonumber
  & \hspace*{6.5em} + \frac{\kappa}{10}\sup_{8r \leq |\alpha|+|\beta| \leq N}
     \int_0^t \left(  \|f\|^2_{\delta-\kappa \tau,l+1,\rho,r,\alpha,\beta}(\tau)
                     +A_{l+1}^2 \right) d\tau \\\nonumber
\leq & \|f(0)\|^2_{\delta,l,\rho,r,\alpha,\beta}
  + C_\kappa \int_0^t \left\{  \|f\|^{2(1+\eta)/\eta}_{l,\rho,r,N}(\tau)
                             + \|f \|^2_{l,\rho,r,N}(\tau)
                      \right\} d\tau
  + C \sup_{\tau \in [0,t]} \sup_{8r\le |\alpha| \le N} \| f \|^2_{\delta-\kappa \tau,l,\rho,r,\alpha+1,0}(\tau)
                       \\\nonumber
  & \hspace*{6.5em} + \frac{\kappa}{10}\sup_{8r \leq |\alpha|+|\beta| \leq N}
     \int_0^t \|f\|^2_{\delta-\kappa \tau,l+1,\rho,r,\alpha,\beta}(\tau)d\tau,
\end{align}
where $0<\eta<1$.

\vspace*{5pt}
\noindent{\bf Step II:}
Now we need to control the term $\| f \|^2_{\delta-\kappa \tau,l,\rho,r,\alpha+1,0}$ with $|\alpha| \le N-1$. For sake of simplicity, we will treat with $\| f \|^2_{\delta-\kappa \tau,l,\rho,r,\alpha,0}\ (|\alpha| \le N)$.

Indeed, it is a matter repeating the argument that we did for the case $|\beta| \ge 1$, except that we need to replace the weight $\frac{\rho^{2(|\alpha|+|\beta|)}}{\{(\alpha-r)!\}^{2\nu_1}\{(\beta-r)!\}^{2\nu_2}}$ by $\frac{\rho^{2|\alpha|}}{\{(\alpha-r)!\}^{2\nu_1}}$, and we can get
\begin{align}\label{Gronw_(alpha,0)}
  &\| f(t)\|^2_{\delta-\kappa t,l,\rho,r,\alpha,0}
 + 2 c_0 \int_0^t
      \frac{\rho^{2|\alpha|} \| W_{l+\frac{\gamma}{2}} \mu^{-1} \partial^\alpha f(\tau)\|^2_{L^2_x(H^s_v)}} {\{ (\alpha-r)!\}^{2\nu_1}} d\tau
 + 2\kappa \int_0^t \| f(\tau) \|^2_{\delta-\kappa \tau,l+1,\rho,r,\alpha,0} d\tau  \\\nonumber
\leq & \|f(0)\|^2_{\delta,l,\rho,r,\alpha,0}
  + C_\kappa \int_0^t \Big\{  \sup_{8r\le |\alpha|\le N} \|f\|^{2(1+\eta)/\eta}_{\delta-\kappa \tau,l+1,\rho,r,\alpha,0}(\tau)
       + A_l^{2/\eta} \sup_{8r\le |\alpha|\le N} \|f\|^2_{\delta-\kappa \tau,l,\rho,r,\alpha,0}(\tau) \\\nonumber
     & \hspace*{12em}
       + \sup_{8r\le |\alpha|\le N} \|f\|^2_{\delta-\kappa \tau,l,\rho,r,\alpha,0}(\tau)
                     \sup_{\substack{8r\le |\alpha|+|\beta| \le N \\ |\beta|=1}} \|f\|^2_{\delta-\kappa \tau,l,\rho,r,\alpha,1}(\tau) \\\nonumber
     & \hspace*{12em}
       + A_l^2 \sup_{\substack{8r\le |\alpha|+|\beta| \le N \\ |\beta|=1}} \|f\|^2_{\delta-\kappa \tau,l,\rho,r,\alpha,1}(\tau)
       + A_l^4
                      \Big\} d\tau \\\nonumber
  & \hspace*{6.5em} + \frac{\kappa}{10}\sup_{8r \leq |\alpha|\leq N}
     \int_0^t \left(  \|f\|^2_{\delta-\kappa \tau,l+1,\rho,r,\alpha,0}(\tau)
                     +A_{l+1}^2 \right) d\tau \\\nonumber
\leq & \|f(0)\|^2_{\delta,l,\rho,r,\alpha,0}
  + C_\kappa \int_0^t \Big\{  \sup_{8r\le |\alpha|\le N} \|f\|^{2(1+\eta)/\eta}_{\delta-\kappa \tau,l+1,\rho,r,\alpha,0}(\tau)
       + \sup_{8r\le |\alpha|\le N} \|f\|^2_{\delta-\kappa \tau,l,\rho,r,\alpha,0}(\tau) \\\nonumber
     & \hspace*{12em}
       + \sup_{8r\le |\alpha|\le N} \|f\|^2_{\delta-\kappa \tau,l,\rho,r,\alpha,0}(\tau)
                     \sup_{\substack{8r\le |\alpha|+|\beta| \le N \\ |\beta|=1}} \|f\|^2_{\delta-\kappa \tau,l,\rho,r,\alpha,1}(\tau) \\\nonumber
     & \hspace*{12em}
       + \sup_{\substack{8r\le |\alpha|+|\beta| \le N \\ |\beta|=1}} \|f\|^2_{\delta-\kappa \tau,l,\rho,r,\alpha,1}(\tau)
                      \Big\} d\tau \\\nonumber
  & \hspace*{6.5em} + \frac{\kappa}{10}\sup_{8r \leq |\alpha|\leq N}
     \int_0^t \|f\|^2_{\delta-\kappa \tau,l+1,\rho,r,\alpha,0}(\tau) d\tau,
\end{align}
where we emphasize especially that, the integral including $\mathcal{J}_2(t)$ in the case of $(\alpha,0)$ can be bounded by
\begin{align}
  & C_\kappa \int_0^t \Big\{
       \sup_{8r\le |\alpha|\le N} \|f\|^2_{\delta-\kappa \tau,l,\rho,r,\alpha,0}(\tau)
           \sup_{\substack{8r\le |\alpha|+|\beta| \le N \\ |\beta|=1}}
           \|f\|^2_{\delta-\kappa \tau,l,\rho,r,\alpha,1}(\tau)
     + A_l^2 \sup_{8r\le |\alpha|\le N} \|f\|^2_{\delta-\kappa \tau,l,\rho,r,\alpha,0}(\tau) \\\nonumber
     & \hspace*{4em}
       + A_l^2 \sup_{\substack{8r\le |\alpha|+|\beta| \le N \\ |\beta|=1}} \|f\|^2_{\delta-\kappa \tau,l,\rho,r,\alpha,1}(\tau)
       + A_l^4
                      \Big\} d\tau \\\nonumber
  & + \frac{\kappa}{20}\sup_{8r \leq |\alpha|\leq N}
     \int_0^t \|f\|^2_{\delta-\kappa \tau,l+1,\rho,r,\alpha,0}(\tau) d\tau.
\end{align}

As a consequence, we have
\begin{align}
  &\| f(t)\|^2_{\delta-\kappa t,l,\rho,r,\alpha,0}
 + 2 c_0 \int_0^t
      \frac{\rho^{2|\alpha|} \| W_{l+\frac{\gamma}{2}} \mu^{-1} \partial^\alpha f(\tau)\|^2_{L^2_x(H^s_v)}} {\{ (\alpha-r)!\}^{2\nu_1}} d\tau \\\nonumber
\leq & \|f(0)\|^2_{\delta,l,\rho,r,\alpha,0}
  + C_\kappa \int_0^t \Big\{  \sup_{8r\le |\alpha|\le N} \|f\|^{2(1+\eta)/\eta}_{\delta-\kappa \tau,l+1,\rho,r,\alpha,0}(\tau)
       + \sup_{8r\le |\alpha|\le N} \|f\|^2_{\delta-\kappa \tau,l,\rho,r,\alpha,0}(\tau) \\\nonumber
     & \hspace*{12em}
       + \sup_{8r\le |\alpha|\le N} \|f\|^2_{\delta-\kappa \tau,l,\rho,r,\alpha,0}(\tau)
                     \sup_{\substack{8r\le |\alpha|+|\beta| \le N \\ |\beta|=1}} \|f\|^2_{\delta-\kappa \tau,l,\rho,r,\alpha,1}(\tau) \\\nonumber
     & \hspace*{12em}
       + \sup_{\substack{8r\le |\alpha|+|\beta| \le N \\ |\beta|=1}} \|f\|^2_{\delta-\kappa \tau,l,\rho,r,\alpha,1}(\tau)
                      \Big\} d\tau,
\end{align}

Thanks to the facts $\{8r\le |\alpha|+|\beta| \le N,\ |\beta|=0\} \subset \{8r\le |\alpha|+|\beta| \le N\}$ and $\{8r\le |\alpha|+|\beta| \le N,\ |\beta|=1\} \subset \{8r\le |\alpha|+|\beta| \le N\}$, we get
\begin{align}\label{Gronw_(alpha,0)_final}
  &\| f(t)\|^2_{\delta-\kappa t,l,\rho,r,\alpha,0}
 + 2 c_0 \int_0^t
      \frac{\rho^{2|\alpha|} \| W_{l+\frac{\gamma}{2}} \mu^{-1} \partial^\alpha f(\tau)\|^2_{L^2_x(H^s_v)}} {\{ (\alpha-r)!\}^{2\nu_1}} d\tau \\\nonumber
\leq & \|f(0)\|^2_{\delta,l,\rho,r,\alpha,0}
  + C \int_0^t \left\{  \|f\|^{2(1+\eta)/\eta}_{l,\rho,r,N}(\tau)
                             + \|f \|^2_{l,\rho,r,N}(\tau)
                      \right\} d\tau.
\end{align}

\vspace*{5pt}
\noindent{\bf Step III:}
Plugging this into \eqref{Gronw} entails that
\begin{align}\label{final inequa}
  &\| f(t)\|^2_{\delta-\kappa t,l,\rho,r,\alpha,\beta}
 + 2 c_0 \int_0^t
      \frac{\rho^{2(|\alpha|+|\beta|)} \| W_{l+\frac{\gamma}{2}} \mu^{-1} \partial^\alpha_\beta f(\tau)\|^2_{L^2_x(H^s_v)}} {\{ (\alpha-r)!\}^{2\nu_1} \{ (\beta-r)!\}^{2\nu_2}} d\tau
 + 2\kappa \int_0^t \| f(\tau) \|^2_{\delta-\kappa \tau,l+1,\rho,r,\alpha,\beta} d\tau  \\\nonumber
\leq & \|f(0)\|^2_{\delta,l,\rho,r,\alpha,\beta}
  + C_\kappa \int_0^t \left\{  \|f\|^{2(1+\eta)/\eta}_{l,\rho,r,N}(\tau)
                             + \|f \|^2_{l,\rho,r,N}(\tau)
                      \right\} d\tau
  + \frac{\kappa}{10}\sup_{8r \leq |\alpha|+|\beta| \leq N}
     \int_0^t \|f\|^2_{\delta-\kappa \tau,l+1,\rho,r,\alpha,\beta}(\tau)d\tau.
\end{align}

This leads to the desired estimate (\ref{basic inequa}) including the extra second term of the left-hand side and so completes the whole proof of Lemma \ref{Key lemma}. \qed

\section{The orders of Gevrey regularity}

We hope to show the orders of Gevrey regularity are $1/s$ for $v$ and $1$ for $x$.

Firstly we modify the definition (\ref{double norm}) as
\begin{align}\label{modified norm}
   \| f \|_{\delta,l,\rho_1,\rho_2,r,\alpha,\beta} \xlongequal{\!\!def\!\!}
   \frac{\rho_1^{|\alpha|}\rho_2^{|\beta|} \| \langle v \rangle^l e^{\delta \langle v \rangle^2}
   \partial^\alpha_\beta f\|_{L^2_x(\mathbb{T}^3) L^2_v(\mathbb{R}^3)} }
   {\{ (\alpha-r)!\}^{\nu_1} \{ (\beta-r)!\}^{\nu_2}}.
\end{align}

\begin{rema}
If $\rho_1=\rho_2=\rho$, the above definition goes back to the definition (\ref{double norm}). Obviously, the previous results are fit for the modified norm, which is just the reason why we only use $\rho$ in the previous process.
\end{rema}

We now introduce a new norm
\begin{align}
|\| f \||_{l,\rho,r,N}(t) \xlongequal{\!\!def\!\!}
\sup_{|\alpha|+|\beta| \leq N} \|f(t)\|_{\delta-\kappa t,l,\rho_1,\rho_2,r,\alpha,\beta},
\end{align}
where $\rho_1=\rho \,(1-t^{\frac{1}{\nu_1}})^{\nu_1}$ and $\rho_2=\rho t^{\nu_2}$.
It should be pointed out that, by the assumption of Theorem \ref{Gevrey Regularity} that $f(0,x,v)$ is analytic in the variable $x$ and the definition of the smooth Maxwellian decay solution (see Def.\ref{Def2}), we have, for sufficiently small $\rho$,
\begin{align}\label{new norm-1-initial}
  |\| f \||_{l,\rho,r,N}(0)
= \sup_{|\alpha|\leq N}
  \frac{\rho^{|\alpha|} \|\langle v \rangle^l e^{\delta \langle v \rangle^2}
        \partial^\alpha_x f(0)\|_{L^2_{x,v}}}
       {\{ (\alpha-r)!\}^{\nu_1}}
= \sup_{|\alpha|\leq N}
  \frac{\rho^{|\alpha|} \|W_l \mu^{-1}(0) \partial^\alpha_x f(0)\|_{L^2_{x,v}}}
       {\{ (\alpha-r)!\}^{\nu_1}} \le C.
\end{align}

As for the above norm, we can get the following integral corresponding to (\ref{EQ2}):
\begin{align}
&\|f(t)\|^2_{\delta -\kappa t,l,\rho_1,\rho_2,r,\alpha,\beta}
+ 2\kappa \!\int_0^t \!
    \|f(\tau)\|^2_{\delta-\kappa \tau,l+1,\rho_1,\rho_2,r,\alpha,\beta}d\tau \\
&\quad + \int_0^t \frac{(2|\alpha|) \rho^{2(|\alpha|+|\beta|)}
   \tau^{\frac{1}{\nu_1}-1} (1-\tau^{\frac{1}{\nu_1}})^{2|\alpha|\nu_1 -1} \tau^{2|\beta|\nu_2}}
                       {\{(\alpha-r)!\}^{2\nu_1}\{(\beta-r)!\}^{2\nu_2}}
                  \| W_l \mu^{-1} \partial^\alpha_\beta f \|^2_{L^2_{x,v}} d\tau \nonumber \\
\leq
&\|f(0)\|^2_{\delta,l,\rho_1,\rho_2,r,\alpha,\beta}
 +2\int^t_0\! \frac{\rho^{2(|\alpha|+|\beta|)} (1-\tau^{\frac{1}{\nu_1}})^{2|\alpha|\nu_1} \tau^{2|\beta|\nu_2}}{\{(\alpha-r)!\}^{2\nu_1}\{(\beta-r)!\}^{2\nu_2}}
   \Big( \mathcal{R}(\tau) + \Psi(\tau) + \mathcal{J}(\tau) + \mathcal{K}(\tau)  \Big) d\tau \nonumber \\
&\quad + \int_0^t \frac{(2|\beta|\nu_2) \rho^{2(|\alpha|+|\beta|)} (1-\tau^{\frac{1}{\nu_1}})^{2|\alpha|\nu_1} \tau^{2|\beta|\nu_2-1}}
                       {\{(\alpha-r)!\}^{2\nu_1}\{(\beta-r)!\}^{2\nu_2}}
                  \| W_l \mu^{-1} \partial^\alpha_\beta f \|^2_{L^2_{x,v}} d\tau,\nonumber
\end{align}
because of the formula
\begin{align*}
&(1-t^{\frac{1}{\nu_1}})^{2|\alpha|\nu_1} \cdot t^{2|\beta|\nu_2} \frac{d}{dt}\|g(t)\|^2 \\
=&\frac{d}{dt} \left( (1-t^{\frac{1}{\nu_1}})^{2|\alpha|\nu_1} \cdot t^{2|\beta|\nu_2} \|g(t)\|^2 \right)
-2|\beta| \nu_2 \cdot (1-t^{\frac{1}{\nu_1}})^{2|\alpha|\nu_1} \cdot t^{2|\beta|\nu_2-1} \|g(t)\|^2 \\ &+ 2|\alpha| \cdot t^{\frac{1}{\nu_1}-1} (1-t^{\frac{1}{\nu_1}})^{2|\alpha|\nu_1 -1} \cdot t^{2|\beta|\nu_2} \|g(t)\|^2.
\end{align*}
Note that, when considering a sufficiently small $\rho$,
\begin{align}
   \|f(0)\|_{\delta,l,\rho_1,\rho_2,r,\alpha,\beta}
 \le |\| f \||_{l,\rho,r,N}(0)
 \le C
\end{align}
is well-defined.

Similar to the argument of Section 4 in \cite{Mori-Ukai}, and noticing the following interpolation inequality
$$
\|f\|^2_{H^k_p(\mathbb{R}^3)} \leq C_\delta
\|f\|_{H^{k-\delta}_{2p}(\mathbb{R}^3)} \|f\|_{H^{k+\delta}_0(\mathbb{R}^3)},
\quad (k \in \mathbb{R},\ p \in \mathbb{R}_+,\ \delta > 0)
$$
implies that for $\gamma \in (-1-2s,\ 1-2s) \subset (-2,\ 1)$,
$$
  \|W_l F^{(\alpha,\beta)} \|^2_{H^{s/2}}
\lesssim  \|W_{l+\gamma/2} F^{(\alpha,\beta)} \|_{H^s}
          \|W_{l-\gamma/2} F^{(\alpha,\beta)} \|_{L^2}
\lesssim  \|W_{l+\gamma/2} F^{(\alpha,\beta)} \|_{H^s}
          \|W_{l+1} F^{(\alpha,\beta)} \|_{L^2},
$$
we can obtain finally $\nu_2=\frac{1}{s}$.

Considering the third term on the left-hand side of the above equation, with $\nu_1=1$ we have
$$
t^{\frac{1}{\nu_1}-1} (1-t^{\frac{1}{\nu_1}})^{2|\alpha|\nu_1 -1}
=(1-t)^{2|\alpha| -1}.
$$
Then the third term can be ignored. Thus, one can obtain immediately the inequality similar to (\ref{final inequa}), which yields the Gevrey smoothing effect in a short interval. That leads us to the conclusion together with Theorem \ref{Propagation of Gevrey}. \qed

\vspace*{2em} \noindent\textbf{Acknowledgements}. This work was
partially supported by NNSFC (No. 11271382 and No. 10971235), RFDP
(No. 20120171110014), and the key project of Sun Yat-sen
University (No. c1185). The authors thank the referee for valuable comments and
suggestions.

\phantomsection

\end{document}